\documentclass[a4paper,11pt]{article}

\usepackage[a4paper]{geometry}
\usepackage{color}
\usepackage{graphicx}
\usepackage{amsmath}
\usepackage{amsfonts}
\usepackage{amsthm}
\usepackage{amscd}
\usepackage[mathcal]{euscript}
\usepackage{epsfig}
\usepackage{amssymb}
\usepackage{slashed}
\usepackage{tabularx}
\usepackage{calligra}
\usepackage{enumerate}
\usepackage{hyperref}
\usepackage{float}
\usepackage{stackrel}
\usepackage[T1]{fontenc}
\usepackage{longtable}
\usepackage{amsthm}
\usepackage{mathrsfs}
\usepackage{verbatim} % for comments
\usepackage{fancyhdr}
\usepackage{tikz}
\usepackage{tikz-cd}
\usepackage{slashed}
\usetikzlibrary{matrix}
\usetikzlibrary{positioning}
\usepackage[all]{xy}
\usepackage{units}
\usepackage{textcomp}
\usepackage{turnstile}
\usepackage{vmargin}
\usepackage{anysize}
\usepackage{ stmaryrd }

\usepackage[font=small]{caption}

\setmargins{2,5cm}       % left margin
{1,5cm}                        %  top margin
{16cm}                      %  text widht
{23,42cm}                    % text height
{10pt}                           % header height
{1cm}                           % space between text and header
{0pt}                             % footer height
{2cm}                           %  space between text and footer

\theoremstyle{plain}

\newtheorem{proposition}{Proposition}

\theoremstyle{definition}

\newcommand{\K}{\mathcal{K}}

\newcommand{\sA}{\mathsf{A}}
\newcommand{\sB}{\mathsf{B}}
\newcommand{\sC}{\mathsf{C}}
\newcommand{\J}{\mathcal{J}}
\renewcommand{\K}{\mathcal{K}}

\begin{document}

\title{On self-dual Yang-Mills fields on special complex surfaces}

\author{Bernardo Araneda\footnote{Email: \texttt{bernardo.araneda@aei.mpg.de}} \\
Max-Planck-Institut f\"ur Gravitationsphysik \\ 
(Albert Einstein Institut), Am M\"uhlenberg 1, \\
D-14476 Potsdam, Germany}

\date{\today}

\maketitle

\begin{abstract}
We derive a generalization of the flat space Yang's and Newman's equations for self-dual Yang-Mills fields 
to (locally) conformally K\"ahler Riemannian 4-manifolds.
The results also apply to Einstein metrics (whose full curvature is not necessarily self-dual).
We analyse the possibility of hidden symmetries in the form of B\"acklund transformations,
and we find a continuous group of hidden symmetries only for the case in which 
the geometry is conformally half-flat. No isometries are assumed.
\end{abstract}

\section{Introduction}

Yang's equation \cite{Yang} is an approach to the analysis of the self-dual Yang-Mills (SDYM) equations in four-dimensional 
flat space, that is intimately related to integrability, solution generating techniques, and twistor theory \cite{MW, MW2}.
It has an infinite dimensional hidden symmetry group, generated by B\"acklund transformations 
combined with conjugation by constant matrices.
The construction is based on the special structure of flat space, so it does not extend to arbitrary curved manifolds.
However, the analysis can be understood as based on the choice of a complex structure, 
so a natural question is whether a similar construction exists for Hermitian manifolds.
The purpose of this note is to address this question,
motivated by the relation of this approach with integrability and solution generating techniques, 
and by recent results \cite{ChenTeo, Tod2020, AA} concerning the relevance of conformally K\"ahler geometry 
to gravitational instantons.

Our main results are as follows:
We find that a Yang's matrix and its field equation can be obtained for conformally K\"ahler manifolds; 
that there is a natural candidate for a B\"acklund transformation; 
and that such a transformation produces a new solution if and only if the geometry is conformally half-flat. 
Combining then with conjugation by constant matrices, we obtain a continuous group of hidden symmetries
for SDYM fields on the conformal class of scalar-flat K\"ahler complex surfaces.

We also obtain a generalization of the closely related `Newman equation' or `$\mathcal{K}$-matrix equation' 
(obtained in \cite{Newman} for flat space) to conformal K\"ahler geometry,
and relate it to other equations in the literature, especially in connection with what is known as the `double copy'.

Following the observation of Atiyah, Hitchin and Singer in \cite{Atiyah} 
that the Einstein property of the metric is equivalent to self-duality of the spin bundles 
(not of the full geometry), we also find Yang and $\K$-matrix equations for the Einstein equations.
In the conformally half-flat case, this allows to generate ${\rm GL}(2,\mathbb{C})$-SDYM fields
using a quaternionic-K\"ahler metric as a seed.
Whether these fields give origin to new quaternionic-K\"ahler metrics is not clear.

Solution generating techniques for different field equations
have been extensively studied in the literature, in a wide variety of contexts.
For curved geometries, these methods frequently assume the existence of isometries.
For example, in the context of the vacuum Einstein equations, a result very similar to the flat space Yang's 
matrix construction exists in general relativity, namely the fact that a reduction of Yang's equation by two commuting symmetries is 
equivalent to the Ernst equation for stationary axisymmetric spacetimes (see \cite{MW2},
\cite[Section 6.6]{MW} and references therein).
`Yang's matrix' in this case is the Gram matrix of Killing fields, and it does not coincide with the Yang's matrix we use in this note. 
Although our curvature restrictions needed for the existence of a hidden symmetry group are rather strong, 
the construction in this work is independent of whether or not there are isometries.

\bigskip
\noindent
{\bf Acknowledgements.} 
I am very grateful to Lionel Mason for comments on this note.
I also gratefully acknowledge the support of the Alexander von Humboldt Foundation.

\section{Background}

The material in this section is not new: it serves to establish our notation and context. 
Our main reference is the book of Mason and Woodhouse \cite{MW}. 
(There are a few places where our conventions differ slightly from those of \cite{MW}, 
but this should not cause any confusion.)
Let $M$ be a 4-dim. smooth manifold, and let $E\to M$ be a rank-$n$ complex vector bundle. 
Unless otherwise specified, we assume the structure group to be ${\rm GL}(n,\mathbb{C})$.
We will be interested in connections on $E$: linear operators ${\rm D}:\Gamma(E)\to\Gamma(T^{*}M\otimes E)$ 
such that ${\rm D}(fs)={\rm d}f\otimes s + f{\rm D}s$ for all smooth functions $f$ and smooth sections $s$.
We also use the notation $({\rm D}s)(X)=X\lrcorner{\rm D}s={\rm D}_{X}s$, for any vector field $X$.
Recall that a local trivialization is a choice of local frame field $\{e_{\bf i}\}$, ${\bf i}=1,...,n$.
Indices ${\bf i}, {\bf j},...$ are concrete; and we also use the Einstein summation convention.
In a local trivialization, we have ${\rm D}e_{\bf i} = A^{\bf j}{}_{\bf i}\otimes e_{\bf j}$ for some matrix-valued 
1-form $A$, called `gauge potential' or `connection 1-form'. 
A local section can be expressed as $s=s^{\bf i}e_{\bf i}$, so 
${\rm D}s = ({\rm d}s^{\bf i} + A^{\bf i}{}_{\bf j}s^{\bf j})\otimes e_{\bf i}$,
which we abbreviate by saying that locally we can write ${\rm D}={\rm d} + A$.

The curvature of ${\rm D}$ is an ${\rm End}(E)$-valued 2-form, $F\in\Gamma(\Lambda^{2}\otimes{\rm End}(E))$, 
defined by $F(X,Y)s={\rm D}_{X}{\rm D}_{Y}s - {\rm D}_{Y}{\rm D}_{X}s - {\rm D}_{[X,Y]}s$, 
where $[\cdot,\cdot]$ is the Lie bracket of vector fields. (We denote by $\Lambda^{k}$ the space of $k$-forms over $M$.)
In a local trivialization, we get the standard expression $F = {\rm d}A + A\wedge A$, 
where $A\wedge A$ is an abbreviation for the matrix-valued 2-form 
defined by $(A\wedge A)^{\bf i}{}_{\bf j}(X,Y)=A(X)^{\bf i}{}_{\bf k}A(Y)^{\bf k}{}_{\bf j} - A(Y)^{\bf i}{}_{\bf k}A(X)^{\bf k}{}_{\bf j}$.
The operator ${\rm D}$ extends naturally to ${\rm End}(E)$-valued differential forms: 
if $\alpha$ is an ${\rm End}(E)$-valued $k$-form, then 
${\rm D}\alpha = {\rm d}\alpha + A\wedge \alpha - (-1)^{k}\alpha\wedge A$.
When applied to the curvature, this gives the Bianchi identity ${\rm D}F=0$.

Let $(z,w,\tilde{z},\tilde{w})$ be arbitrary local coordinates (later these will be double null coordinates, but for 
the moment they are arbitrary ---notice that we have not yet introduced a metric). 
The associated coordinate basis for the tangent space is 
$(\partial_{z},\partial_{w},\partial_{\tilde{z}},\partial_{\tilde{w}})$, and the cotangent basis is 
$({\rm d}z, {\rm d}w, {\rm d}\tilde{z}, {\rm d}\tilde{w})$. 
The exterior derivative can be decomposed as ${\rm d} = \partial +\tilde{\partial}$, where $\partial={\rm d}z\partial_{z}+{\rm d}w\partial_{w}$ 
and $\tilde{\partial}={\rm d}\tilde{z}\partial_{\tilde{z}}+{\rm d}\tilde{w}\partial_{\tilde{w}}$.
We can write $A = A_{z}{\rm d}z + A_{w}{\rm d}w + A_{\tilde{z}}{\rm d}\tilde{z} + A_{\tilde{w}}{\rm d}\tilde{w}$, 
where $A_{z} = \partial_{z} \lrcorner A$,  $A_{w} = \partial_{w} \lrcorner A$, etc. are matrix-valued functions.
Similarly, $F = F_{zw}{\rm d}z\wedge{\rm d}w + ...$, where $F_{zw} = F(\partial_{z},\partial_{w})$,... are matrices.
In terms of $A$, we have for example
\begin{align}
 F_{zw} ={}& \partial_{z}A_{w} - \partial_{w}A_{z} + [A_{z},A_{w}], \label{Fzw} \\
 F_{\tilde{z}\tilde{w}} ={}& \partial_{\tilde{z}}A_{\tilde{w}} - \partial_{\tilde{w}}A_{\tilde{z}} + [A_{\tilde{z}},A_{\tilde{w}}], \label{Ftztw}
\end{align}
and similarly for the other components, where now $[\cdot,\cdot]$ denotes the commutator of matrices.

Assume now that the manifold is equipped with a metric and an orientation. 
We then have a Hodge star operator $*$, which maps 2-forms to 2-forms. 
In Riemann signature it holds $*^{2}=1$, so $*$ has eigenvalues $\pm 1$ and 
produces a decomposition $\Lambda^{2}=\Lambda^{2}_{+}\oplus\Lambda^{2}_{-}$ into 
self-dual (SD) and anti-self-dual (ASD) spaces. 
The Yang-Mills equations are ${\rm D}*F=0$, and the self-dual Yang-Mills equations are
\begin{equation}
 *F = F. \label{SDYMeq}
\end{equation}
These equations are gauge invariant, so one can analyse them in any gauge.
A convenient approach to the study of \eqref{SDYMeq} is to find a basis of ASD 2-forms, 
and then impose the corresponding ASD components of $F$ to vanish. 
Although the discussion of the previous paragraphs, together with expressions like \eqref{Fzw}, \eqref{Ftztw}, 
do not depend on the background metric structure, which 2-forms are ASD does of course depend on it.

\bigskip
For example, it is instructive to briefly look at the simple case of complexified Minkowski spacetime $\mathbb{CM}$.
Assume that $(z,w,\tilde{z},\tilde{w})$ are double null coordinates: 
the metric is $\eta=2({\rm d}z \: {\rm d}\tilde{z} - {\rm d}w \: {\rm d}\tilde{w})$. 
A basis of ASD 2-forms\footnote{Our conventions for SD and ASD forms is opposite to the one in \cite{MW}.} 
is ${\rm d}z\wedge{\rm d}w$, ${\rm d}\tilde{z}\wedge{\rm d}\tilde{w}$, 
${\rm d}z\wedge{\rm d}\tilde{z} - {\rm d}w\wedge{\rm d}\tilde{w}$.
Therefore, the SDYM equations \eqref{SDYMeq} have in this case the well-known form
\begin{equation}
 F_{zw} = 0 , \qquad F_{\tilde{z}\tilde{w}} = 0, \qquad F_{z\tilde{z}} - F_{w\tilde{w}} = 0. \label{SDYMflatspace}
\end{equation}
In terms of the gauge potential, the first two equations are given by equating to zero the right-hand-sides of 
\eqref{Fzw}, \eqref{Ftztw}. 
As explained in \cite{MW}, these are the local integrability conditions for the existence of 
two frames $\{e_{\bf i}\}$ and $\{\tilde{e}_{\bf i}\}$ such that 
\begin{align}
 & {\rm D}_{z} e_{\bf i} = 0, \qquad {\rm D}_{w} e_{\bf i} = 0, \label{specialframe} \\
 & {\rm D}_{\tilde{z}}\tilde{e}_{\bf i} = 0, \qquad {\rm D}_{\tilde{w}}\tilde{e}_{\bf i} = 0 \label{specialframetilde}
\end{align}
for all ${\bf i}$, where ${\rm D}_{z} = \partial_{z}\lrcorner {\rm D}$, etc.
The corresponding gauge potentials, $A$ and $\tilde{A}$, satisfy $A_{z}=0=A_{w}$ and $\tilde{A}_{\tilde{z}}=0=\tilde{A}_{\tilde{w}}$.
Define now a matrix $\J=( \J^{\bf i}{}_{\bf j})$ by $e_{\bf j} = \J^{\bf i}{}_{\bf j}\tilde{e}_{\bf i}$. This matrix is called {\em Yang's matrix}.
In terms of $\J$, we have $A = \J^{-1}\tilde{\partial}\J$ and $\tilde{A}=\J\partial\J^{-1}$.
The first two equations in \eqref{SDYMflatspace} are identically satisfied since 
$A_{z}=0=A_{w}$, $\tilde{A}_{\tilde{z}}=0=\tilde{A}_{\tilde{w}}$, and the equations \eqref{SDYMeq} are gauge invariant.
The third equation is satisfied if and only if $\J$ satisfies {\em Yang's equation}
\begin{equation}
 \partial_{z}(\J^{-1}\partial_{\tilde{z}}\J) - \partial_{w}(\J^{-1}\partial_{\tilde{w}}\J)  = 0. \label{YangEqFlatspace}
\end{equation}

This equation has `hidden symmetries' in the form of B\"acklund transformations \cite{MW}, 
which produce new solutions to the SDYM equations \eqref{SDYMflatspace}. 
The solutions are gauge inequivalent, see \cite{MCN}. (Nevertheless, they can be understood in terms of 
`irregular gauge transformations' that depend on a spectral parameter 
$\zeta\in\mathbb{CP}^{1}=\mathbb{C}\cup\{\infty\}$ 
---while a standard gauge transformation does not---
and are not well-behaved at $\zeta=0$ and $\zeta=\infty$, see \cite[Section 4.6]{MW}.)
By combining B\"acklund transformations with conjugation of $\J$ by constant matrices, 
an infinite-dimensional group of symmetries is generated.
The B\"acklund transformation has a simple description in twistor space $\mathbb{PT}$, 
as a special conjugation of the patching matrix that defines a holomorphic vector bundle over $\mathbb{PT}$ 
associated to the SDYM field via the Penrose-Ward transform; see \cite[Section 12.1]{MW}.

The Yang's matrix approach in flat space, together with its connection with twistor theory, 
can be conveniently formulated in terms of the Lax pair 
$L={\rm D}_{\tilde{z}}-\zeta{\rm D}_{w}$, $M={\rm D}_{\tilde{w}}-\zeta{\rm D}_{z}$, 
where $\zeta\in\mathbb{CP}^{1}$ is the spectral parameter mentioned before:
the SDYM equation \eqref{SDYMeq} is equivalent to the condition $[L,M]=0$ for all $\zeta$.
This formulation does not seem to be convenient for the purposes of this note, 
since we want to consider the case where there is only one integrable complex structure, 
which is equivalent to saying that only one value of the spectral parameter is special\footnote{Nevertheless, 
we point out that a formulation in terms of a different Lax pair and a spectral parameter might still be possible, 
as is done in \cite{BV} for the case of the Einstein vacuum equations with two commuting symmetries.}.
More precisely as follows. Recall that the space of (orthogonal, almost-) complex 
structures in four dimensions in Riemann signature is the projective spin bundle $\mathbb{PS}$, whose
fibers are complex projective lines $\mathbb{CP}^{1}$. 
The spectral parameter can then be understood as a coordinate on the fibers of $\mathbb{PS}$, 
and it can be equivalently represented by a projective spinor, say $[\zeta^{A}]$. 
(We use abstract spacetime and spinor indices following \cite{PR1}.)
The complex structure $J_{\zeta}$ defined by a particular $[\zeta^{A}]$ is 
\begin{equation}
 (J_{\zeta})^{a}{}_{b} = \frac{i}{(\zeta_{C}\zeta^{\dagger C})}(\zeta^{A}\zeta^{\dagger}_{B}+\zeta^{\dagger A}\zeta_{B})\delta^{A'}_{B'}
 \label{Jzeta}
\end{equation}
where $\zeta^{\dagger A}$ is the complex conjugate of $\zeta^{A}$. 
(The symbol $\dagger$ denotes Euclidean spinor conjugation: 
if $\psi^{A}$ has components $(a,b)$, then $\psi^{\dagger A}$ has components $(-\bar{b},\bar{a})$; 
and it holds $\dagger^{2}=-1$.)

Integrability of \eqref{Jzeta} is equivalent to the existence of complex coordinates $(\lambda,\mu)$ such that 
\begin{equation}
 {\rm d}\lambda=\zeta_{A}\alpha_{A'}{\rm d}x^{AA'}, \qquad {\rm d}\mu=\zeta_{A}\beta_{A'}{\rm d}x^{AA'} \label{Coordinates1forms}
\end{equation}
where $\alpha_{A'},\beta_{A'}$ are two spinors with $\alpha_{A'}\beta^{A'}\neq0$. 
We denote the complex conjugate coordinates by $\tilde{\lambda},\tilde{\mu}$.
In (conformally) flat space one can set $\beta_{A'} \propto \alpha^{\dagger}_{A'}$.
The complex 2-surfaces defined by $\{ \lambda={\rm const.}, \; \mu={\rm const.}\}$ 
and $\{ \tilde{\lambda}={\rm const.}, \; \tilde{\mu}={\rm const.}\}$ 
are totally null and are called $\beta$-surfaces: $J_{\zeta}$ defines two 2-parameter families of them.
We call each family a `2-dimensional twistor space'.
Now choose an arbitrary spinor $o^{A}$, let $\iota^{A}=o^{\dagger A}$, and put $\zeta^{A}=o^{A}+\zeta \iota^{A}$.
Replacing in \eqref{Jzeta}, a short calculation gives
\begin{equation}
 J_{\zeta} = \frac{(1-\zeta\bar\zeta)}{(1+\zeta\bar\zeta)}J_{0} + \frac{i(\zeta-\bar\zeta)}{(1+\zeta\bar\zeta)}J_{1} 
 + \frac{(\zeta+\bar\zeta)}{(1+\zeta\bar\zeta)}J_{2} \label{Jzeta2}
\end{equation}
where the triple $(J_{0},J_{1},J_{2})$ satisfies the quaternion algebra. For a fixed value of $\zeta$, the associated complex 
structure is \eqref{Jzeta2}, and its holomorphic coordinates are $(\lambda,\mu)$.

Define now new coordinates by ${\rm d}z = o_{A}\alpha_{A'}{\rm d}x^{AA'}$, ${\rm d}w = o_{A}\beta_{A'}{\rm d}x^{AA'}$,
$\tilde{z}=\bar{z}$, $\tilde{w}=\bar{w}$. The coordinates $(z,w)$ are holomorphic w.r.t the complex structure 
corresponding to $\zeta=0$ in \eqref{Jzeta2}. We have 
\begin{equation}
 \lambda = z + \zeta\tilde{w}, \qquad \mu = w + \zeta\tilde{z}. \label{GeneralHolomorphicCoordinates}
\end{equation}
This establishes the relationship between the coordinates associated to $J_{\zeta}$ 
and the ones associated to the complex structure defined by the spinor $o_{A}$, i.e. $J_{0}$.
The three complex coordinates $(\lambda,\mu,\zeta)$ define a 3-dimensional complex manifold, 
the twistor space $\mathbb{PT}$ of (conformally) flat space.
Twistor space is fibered over $\mathbb{CP}^{1}$: $\zeta$ is a coordinate on the base, 
and $(\lambda,\mu)$ are coordinates on the fiber over $\zeta$.
Each fixed valued of $\zeta$ defines an integrable complex structure $J_{\zeta}$ from the point of view of flat space, 
and two fibers (two 2-dim. twistor spaces) from the point of view of the fibration $\mathbb{PT}\to\mathbb{CP}^{1}$.
The holomorphic coordinates \eqref{GeneralHolomorphicCoordinates} of $J_{\zeta}$ 
parametrize the different $\beta$-surfaces associated to $J_{\zeta}$ on (the complexification of) flat space.

For a generic Hermitian manifold, equations \eqref{Jzeta}, \eqref{Jzeta2} are still valid (since they simply describe 
{\em almost-}complex structures), but there is only one {\em integrable} complex structure:
only one special value of $\zeta$. We can take it to be $\zeta=0$, so $J_{\zeta}=J_{0}$.
The only surviving twistor structures are the two 2-dimensional twistor spaces defined by $J_{0}$.

\section{Conformally K\"ahler manifolds}\label{Sec:curvedspaces}

We will now investigate the possible generalization of the flat space results to conformally K\"ahler Riemannian manifolds.
Our construction will be purely local.

Let $(M,g)$ be a 4-dimensional orientable Riemannian manifold. 
Assume that there is (locally) an integrable, orthogonal almost-complex structure $J$; 
we can take it to be $J_{0}$ in \eqref{Jzeta2}.
(The local condition is necessary since $J$ may not exist globally, e.g. if the manifold is the 4-sphere $S^{4}$.)
Let $\kappa(\cdot,\cdot)=g(J\cdot,\cdot)$ be the fundamental 2-form. One can show that this is SD or ASD.
Since we want to study the SDYM equations, we choose $\kappa$ to be ASD $*\kappa=-\kappa$.
Let $(z,w)$ be local holomorphic coordinates w.r.t. $J$, and $(\tilde{z},\tilde{w})$ local anti-holomorphic coordinates.
This means that the metric, complex structure, and fundamental 2-form are
\begin{align}
 g ={}& 2( g_{z\tilde{z}}{\rm d}z \: {\rm d}\tilde{z} + g_{z\tilde{w}}{\rm d}z \: {\rm d}\tilde{w} 
 + g_{w\tilde{z}}{\rm d}w \: {\rm d}\tilde{z} + g_{w\tilde{w}}{\rm d}w \: {\rm d}\tilde{w}), \label{Hermitianmetric} \\
 J={}& i(\partial_{z}\otimes{\rm d}z + \partial_{w}\otimes{\rm d}w 
 - \partial_{\tilde{z}}\otimes{\rm d}\tilde{z} - \partial_{\tilde{w}}\otimes{\rm d}\tilde{w}), \label{complexstructure} \\
 \kappa ={}& i( g_{z\tilde{z}}{\rm d}z \wedge {\rm d}\tilde{z} + g_{z\tilde{w}}{\rm d}z \wedge {\rm d}\tilde{w} 
 + g_{w\tilde{z}}{\rm d}w \wedge {\rm d}\tilde{z} + g_{w\tilde{w}}{\rm d}w \wedge {\rm d}\tilde{w}). \label{Fundamental2form}
\end{align}
In order to analyse the SDYM equations \eqref{SDYMeq}, we need to find a basis of ASD 2-forms. 
One of them is $\kappa$. 
For the others, we observe that the two eigenbundles of $J$ are totally null (i.e. $g$ vanishes upon restriction to them), 
so the two 2-surface elements must be SD or ASD; see \cite[Lemma 2.3.3]{MW}.
For our choice of orientation (i.e. $*\kappa=-\kappa$), one can show that 
a 2-form $\Omega$ is SD if and only if $\Omega(J\cdot, J\cdot)=\Omega(\cdot,\cdot)$. 
If $\Sigma$ is the surface element of an eigenbundle of $J$ we find $\Sigma(J\cdot, J\cdot)=-\Sigma(\cdot,\cdot)$; 
so it must be ASD (since it is not SD).
In other words, the 2-forms ${\rm d}z\wedge{\rm d}w$ and ${\rm d}\tilde{z}\wedge{\rm d}\tilde{w}$ are ASD. 
(Alternatively, a straightforward proof of this can be done using spinors.)
Therefore, denoting by $F_{\kappa}$ the component of $F$ corresponding to $\kappa$, 
the SDYM equations \eqref{SDYMeq} are now
\begin{equation}
 F_{zw} = 0 , \qquad F_{\tilde{z}\tilde{w}} = 0, \qquad F_{\kappa} = 0. \label{SDYMcurvedspace}
\end{equation}
So we see that the first two equations are exactly the same as in the flat case \eqref{SDYMflatspace}.
We can interpret them as the local integrability conditions for the existence of 
frames $\{e_{\bf i}\}$ and $\{\tilde{e}_{\bf i}\}$ such that equations \eqref{specialframe} and \eqref{specialframetilde} hold.
This leads to the existence of a Yang's matrix $\J=(\J^{\bf i}{}_{\bf j})$ as before, defined by $e_{\bf j} = \J^{\bf i}{}_{\bf j}\tilde{e}_{\bf i}$.
The gauge potential is again given by 
\begin{equation}
 A = \J^{-1}\tilde{\partial}\J. \label{gaugePotential}
\end{equation}
So the generalization of the flat space Yang's equation is $F_{\kappa} = 0$ written in terms of $\J$.
Since now the 2-form ${\rm d}z\wedge{\rm d}\tilde{z} - {\rm d}w\wedge{\rm d}\tilde{w}$ is not ASD, 
this is not simply $F_{z\tilde{z}} - F_{w\tilde{w}} = 0$ as in the flat case.
In fact, using \eqref{Hermitianmetric}--\eqref{Fundamental2form}, one can show that $F_{\kappa}=0$ is equivalent to
\begin{equation}
 g_{w\tilde{w}}F_{z\tilde{z}} + g_{z\tilde{z}}F_{w\tilde{w}} - g_{w\tilde{z}}F_{z\tilde{w}} - g_{z\tilde{w}}F_{w\tilde{z}} = 0.
 \label{F_k}
\end{equation}
This is considerably more difficult than the flat space version.
To find a convenient expression that allows us to analyse hidden symmetries, we will use spinors.
Since we work locally, the existence of a global spin structure is not necessary.

\smallskip
In order to have a geometric interpretation of the spinor operators that will appear below (cf. \eqref{partialconnections}),
it may be instructive to first notice that the equations $F_{zw} = 0$, $F_{\tilde{z}\tilde{w}} = 0$
can be understood using holomorphic vector bundles and partial connections. 
The latter are `connections' that differentiate sections of a vector bundle along vectors lying on a subbundle of $TM$.
First, notice that the complex structure $J$ of the base manifold induces a decomposition of the cotangent bundle 
as $T^{*}M\otimes \mathbb{C}=\Lambda_{(1,0)}\oplus\Lambda_{(0,1)}$, where $\Lambda_{(1,0)}$ (resp. $\Lambda_{(0,1)}$)
corresponds to 1-forms with eigenvalue $+i$ (resp. $-i$) under $J$.
Given a smooth, complex vector bundle $E\to M$, we can define (see \cite[Section 9.5]{MW})
a ``type $(0,1)$'' partial connection as a linear map $\tilde{\partial}_{E}:\Gamma(E)\to\Gamma(\Lambda_{(0,1)}\otimes E)$ 
such that $\tilde{\partial}_{E}(fs)=\tilde{\partial}f\otimes s + f\tilde{\partial}_{E}s$, and similarly 
for a ``type $(1,0)$'' partial connection $\partial_{E}:\Gamma(E)\to\Gamma(\Lambda_{(1,0)}\otimes E)$ 
\footnote{We add the terms ``type $(0,1)$'' and ``type $(1,0)$'' 
to `partial connection' simply to distinguish whether they map to $\Lambda_{(0,1)}$ or to $\Lambda_{(1,0)}$;
in \cite{MW} only the first of these is considered. An alternative terminology, see \cite[Chapter 9]{Moroianu}, 
is to say that $\tilde{\partial}_{E}$ equips $E$ with a `pseudo-holomorphic' structure, 
which will be holomorphic iff $\tilde{\partial}^{2}_{E}=0$.}.
The Yang-Mills connection ${\rm D}$ and the background complex structure $J$ 
give naturally partial connections $\tilde{\partial}_{E}$ and $\partial_{E}$: 
in a local smooth trivialization, using ${\rm d}=\partial+\tilde{\partial}$ and $A = a + \tilde{a}$, 
where $a=A_{z}{\rm d}z+A_{w}{\rm d}w$ is type $(1,0)$
and $\tilde{a}=A_{\tilde{z}}{\rm d}\tilde{z} + A_{\tilde{w}}{\rm d}\tilde{w}$ is type $(0,1)$, 
we can write ${\rm D}=\partial_{E}+\tilde{\partial}_{E}$, with $\partial_{E}=\partial + a$ and 
$\tilde{\partial}_{E}=\tilde{\partial} + \tilde{a}$. 
The vector bundle $E$ will be holomorphic w.r.t. $\tilde{\partial}_{E}$ 
if and only if the following {\em partial flatness} condition holds: $\tilde{\partial}\tilde{a}+\tilde{a}\wedge\tilde{a}=0$.
This condition guarantees that there are local holomorphic trivializations, that is, frames such 
that $\tilde{\partial}_{E}\tilde{e}_{\bf i}=0$.
A short calculation gives $\tilde{\partial}\tilde{a}+\tilde{a}\wedge\tilde{a}=F_{\tilde{z}\tilde{w}}{\rm d}\tilde{z}\wedge{\rm d}\tilde{w}$;
so we can interpret the SDYM equation $F_{\tilde{z}\tilde{w}}=0$ as the condition for the existence of local holomorphic trivializations. 
Similarly, we get $\partial a + a\wedge a = F_{zw}{\rm d}z\wedge{\rm d}w$, 
so $F_{zw}=0$ gives local anti-holomorphic trivializations\footnote{As noted in \cite{MW}, in the general case one needs 
additional data to define the operator $\partial_{E}$; in our case this is provided by the background complex structure and 
the Yang-Mills connection.}.

\subsubsection*{Spinor approach}

We want to investigate whether the third SDYM equation in \eqref{SDYMcurvedspace}, $F_{\kappa}=0$ 
(or more explicitly \eqref{F_k}), 
has a structure similar to the flat space Yang's equation \eqref{YangEqFlatspace}, so that in particular 
we can analyse the possibility of hidden symmetries in the form of B\"acklund transformations.
We will use abstract spacetime and spinor indices, following the conventions of \cite{PR1}.
Denote the spin bundles by $\mathbb{S}$ and $\mathbb{S}'$, equipped with symplectic forms $\epsilon$ and $\epsilon'$.
The isomorphisms $\Lambda^{2}_{+}=\mathbb{S}'^{*}\odot\mathbb{S}'^{*}$ and $\Lambda^{2}_{-}=\mathbb{S}^{*}\odot\mathbb{S}^{*}$
imply that any (possibly matrix-valued) 2-form $F$ can be expressed as $F_{ab}=\psi_{A'B'}\epsilon_{AB}+\varphi_{AB}\epsilon_{A'B'}$, 
where $\psi_{A'B'}=\psi_{(A'B')}$ and $\varphi_{AB}=\varphi_{(AB)}$ are (matrix-valued) symmetric spinors representing 
the SD and ASD parts of $F_{ab}$, respectively. 
For the Yang-Mills curvature $F$, in a local trivialization we have
\begin{align}
 \psi_{A'B'} =& \nabla_{A(A'}A^{A}{}_{B')} + A_{A(A'}A^{A}{}_{B')}, \label{SDYMfield} \\
 \varphi_{AB} =& \nabla_{A'(A}A^{A'}{}_{B)} + A_{A'(A}A^{A'}{}_{B)}, \label{ASDYMfield}
\end{align}
where $\nabla_{AA'}$ is the Levi-Civita connection of $g_{ab}$, and $A_{AA'}$ is the spinor equivalent of 
the gauge potential $A_{a}$.
In this notation, the SDYM equation \eqref{SDYMeq} is $\varphi_{AB} = 0$.

Recall now that the spinor representation of a complex structure is \eqref{Jzeta}.
We work with the specific value $\zeta=0$, i.e. with the spinor $o_{A}$ and its complex conjugate $\iota_{A}=o^{\dagger}_{A}$:
\begin{equation}
 J^{a}{}_{b} = i\chi^{-1}(o^{A}\iota_{B}+\iota^{A}o_{B})\delta^{A'}_{B'}, \qquad \chi:=o_{A}\iota^{A}. \label{ComplexStructureSpinors}
\end{equation}
The fundamental 2-form is $\kappa_{ab}=2i\chi^{-1}o_{(A}\iota_{B)}\epsilon_{A'B'}$. 
The third equation in \eqref{SDYMcurvedspace} is then simply $o^{A}\iota^{B}\varphi_{AB} = 0$. 

The integrability condition on $J^{a}{}_{b}$ is equivalent to the equation $o^{A}o^{B}\nabla_{AA'}o_{B}=0$ on $[o^{A}]$.
The analogue of \eqref{Coordinates1forms} in terms of vector fields is
\begin{align}
 & \partial_{z} = \iota^{A}Z^{A'}\partial_{AA'}, \qquad \partial_{\tilde{z}} = o^{A}\tilde{Z}^{A'}\partial_{AA'}, 
 \qquad Z_{A'}W^{A'}=N, \label{CoordinateVF1} \\
 & \partial_{w} = \iota^{A}W^{A'}\partial_{AA'}, \qquad \partial_{\tilde{w}} = o^{A}\tilde{W}^{A'}\partial_{AA'}, 
 \qquad \tilde{Z}_{A'}\tilde{W}^{A'}=\tilde{N}, \label{CoordinateVF2}
\end{align}
for some spinor fields $(Z^{A'},W^{A'})$ and $(\tilde{Z}^{A'},\tilde{W}^{A'})$ and
some scalar fields $N,\tilde{N}$, where $\partial_{AA'}=\partial/\partial x^{AA'}$.
The determinant of the metric is $\det(g_{\mu\tilde{\nu}})=\chi^{2}N\tilde{N}$.
For the reality conditions $\tilde{z}=\bar{z}$ and $\tilde{w}=\bar{w}$, we get
$\tilde{Z}^{A'}=-Z^{\dagger A'}$ and $\tilde{W}^{A'}=-W^{\dagger A'}$, which give $\tilde{N}=\bar{N}$.
Notice that $g_{z\tilde{w}}=g(\partial_{z},\partial_{\tilde{w}})=\chi Z_{A'}W^{\dagger A'}$, 
which means that in the general case (i.e. for $g_{z\tilde{w}}\neq0$) we have $W^{A'}\neq Z^{\dagger A'}$.
In fact, if one imposes $W^{A'}=\alpha Z^{\dagger A'}$ for some scalar field $\alpha$, then the metric is
$g=g_{z\tilde{z}}({\rm d}z\odot{\rm d}\tilde{z} - |\alpha|^{2}{\rm d}w\odot{\rm d}\tilde{w})$; 
in particular, $\alpha={\rm constant}$ implies conformal flatness.

The equations $F_{zw} = 0$, $F_{\tilde{z}\tilde{w}} = 0$ are simply $\iota^{A}\iota^{B}\varphi_{AB}=0$ 
and $o^{A}o^{B}\varphi_{AB}=0$, respectively, and, as mentioned, they imply the existence of local holomorphic 
and anti-holomorphic trivializations of $E$.

We now impose the (locally) conformally K\"ahler condition: we assume that there is a scalar field $\phi$ such that
\begin{equation}
 {\rm d}(\phi^{2}\kappa) = 0. \label{conformalKahler}
\end{equation}
This will be automatically satisfied if the manifold is Einstein, but here we are not assuming the Einstein condition.
From \cite{DunajskiTod}, the conformal K\"ahler condition implies that the Weyl tensor is one-sided type D.
For practical calculations, it is convenient to express 
equation \eqref{conformalKahler} in terms of the Lee form $f$ of $J$ 
(which is defined by ${\rm d}\kappa = -2f\wedge\kappa$): $f_{a} = \partial_{a}\log\phi$.

To simplify the computations, we choose the scaling of $o^{A}$ so that the following condition holds:
\begin{equation}
 o^{A}\nabla_{AA'}o_{B} = 0. \label{SpecialScaling}
\end{equation}
The compatibility condition for this is $\Psi_{ABCD}o^{B}o^{C}o^{D}=0$ (where $\Psi_{ABCD}$ is the ASD Weyl spinor),
which is satisfied since $\Psi_{ABCD}$ must be type D.
Then a calculation shows that, in terms of the Yang's matrix $\J$ defined above equation \eqref{gaugePotential}, 
the last SDYM equation, $F_{\kappa}=0$, is equivalent to
\begin{equation}
 \nabla^{A'}(r\J^{-1}\tilde{\nabla}_{A'}\J) = 0, \label{YangEqCurvedSpace}
\end{equation}
where for convenience we defined 
\begin{equation}
 r:=\chi^{-1}\phi \label{r}
\end{equation}
together with 
\begin{equation}
 \nabla_{A'}:=\iota^{A}\nabla_{AA'}, \qquad \tilde{\nabla}_{A'}:=o^{A}\nabla_{AA'}. \label{partialconnections}
\end{equation}
These operators correspond to differentiation along type $(1,0)$ and $(0,1)$ vectors:
they are partial connections as described previously in abstract terms.
We summarize the above result in the form of the following:

\begin{proposition}\label{Prop:YangEquation}
Let $(M,g,J)$ be a conformally K\"ahler manifold.
Choose the scaling of the spinor field $o^{A}$ representing $J$ so as to satisfy \eqref{SpecialScaling}, 
and let $\iota_{A}=o^{\dagger}_{A}$.
Let $\J$ be a matrix-valued function, and define a matrix-valued 1-form $A$ by eq. \eqref{gaugePotential}, or in spinor terms:
\begin{equation}
 A_{AA'}=-\chi^{-1}\J^{-1}\iota_{A}o^{B}\nabla_{BA'}\J. \label{gaugePotentialSpinors}
\end{equation}
Then $A_{AA'}$ is the gauge potential of a SDYM field if and only if $\J$ satisfies equation \eqref{YangEqCurvedSpace}.
Conversely, any SDYM field can be represented by a gauge potential of the form \eqref{gaugePotentialSpinors} 
for some matrix $\J$ subject to \eqref{YangEqCurvedSpace}.
\end{proposition}

Equation \eqref{YangEqCurvedSpace} is the generalization of the flat space Yang's equation \eqref{YangEqFlatspace}
to conformal K\"ahler manifolds. In double null coordinates, a calculation shows that \eqref{YangEqCurvedSpace} is
\begin{align}
 \partial_{z}[sg_{w\tilde{z}}\J^{-1}\partial_{\tilde{w}}\J] - \partial_{z}[sg_{w\tilde{w}}\J^{-1}\partial_{\tilde{z}}\J] 
 - \partial_{w}[sg_{z\tilde{z}}\J^{-1}\partial_{\tilde{w}}\J] + \partial_{w}[sg_{z\tilde{w}}\J^{-1}\partial_{\tilde{z}}\J] = 0
 \label{YangEqCurvedSpaceCoordinates}
\end{align}
where for simplicity we put $s=r\tilde{N}^{-1}\chi^{-1}$. For flat space this reduces to \eqref{YangEqFlatspace}.
For the analysis of the possibility of hidden symmetries in the form of B\"acklund transformations,
we find equation \eqref{YangEqCurvedSpace} to be much nicer than the coordinate expression 
\eqref{YangEqCurvedSpaceCoordinates}, so below we will focus on \eqref{YangEqCurvedSpace}.

For the simplest case of Maxwell fields, $\J$ is a scalar field, and \eqref{YangEqCurvedSpace} can be shown to be
the so-called Fackerell-Ipser equation for the field $\Phi=\log\J$. 
(This appears in the context of spin-weight zero electromagnetic perturbations of the Kerr solution in the Lorentzian case.)

\subsubsection*{The $\K$-matrix}

An alternative to the Yang's matrix formulation in the flat case is the so-called `$\K$-matrix' 
and its corresponding equation, see \cite{Newman} and \cite[Section 3.3]{MW}. 
Here we show that this can also be derived in our current context. 
Although we will not explore much about this equation in this note, 
it is interesting to point out that the $\K$-matrix turns out to be connected to other developments 
of past and modern interest, see below.

\smallskip
Start from the SDYM equation $F_{zw}=0$, which gives a local anti-holomorphic trivialization: 
a frame such that $\partial_{E}e_{\bf i}=0$. (Recall that `$E$' here is not a spinor index; it means that
$\partial_{E}$ is a partial connection on the bundle $E$.)
In spinor terms, the associated gauge potential satisfies $\iota^{A}A_{AA'}=0$. 
Now focus on the equation $F_{\kappa}=0$, i.e. $\chi^{-1}o^{A}\iota^{B}\varphi_{AB}=0$.
The same calculation that leads to \eqref{YangEqCurvedSpace} gives $\nabla^{A'}(ro^{A}A_{AA'})=0$.
This implies that there is (locally) a matrix-valued function $\K$ such that $ro^{A}A_{AA'} = \nabla_{A'}\K$.
(The argument for this integration procedure is given below, see the paragraph above eqs. \eqref{B'}-\eqref{tB'}.)
Now replace this in the last SDYM equation: $F_{\tilde{z}\tilde{w}}=0$ (or $o^{A}o^{B}\varphi_{AB}=0$).
Then a calculation gives
\begin{equation}
 \tilde{\nabla}_{A'}(r^{-1}\nabla^{A'}\K) + r^{-2}(\nabla_{A'}\K)(\nabla^{A'}\K) = 0. \label{Kmatrixeq}
\end{equation}
This is a generalization of Newman's equation \cite[Eq. (5.11)]{Newman} 
(or the $\K$-matrix equation \cite[Eq. (3.3.7)]{MW}) to conformally K\"ahler manifolds.
The gauge potential is $A_{AA'}=-\phi^{-1}\iota_{A}\nabla_{A'}\K$.

\smallskip
It turns out that particular cases of \eqref{Kmatrixeq} appear in several a priori unrelated contexts.
For Maxwell fields (or Abelian gauge theories), the term quadratic in $\K$ in \eqref{Kmatrixeq} vanishes, 
and the resulting equation can be shown to be the Teukolsky equation for electromagnetic perturbations 
of type D spaces in general relativity.
In the context of scattering ampitudes in QCD, the flat space version of 
\eqref{Kmatrixeq} appears, in coordinate form, in \cite{Bardeen} (see page 5 in this reference).
For the SD Einstein vacuum equations,
a calculation shows that \eqref{Kmatrixeq} gives Pleba\'nski's second heavenly equation; 
we show this in the appendix \ref{App:2HE}. 
It also appears that \eqref{Kmatrixeq} should reduce to the hyper-heavenly equation of Pleba\'nski and 
Robinson, or to the `conformal HH equation' for conformal geometry derived in \cite{Araneda}.

It is interesting to note that Pleba\'nski's heavenly equation can be obtained 
from \eqref{Kmatrixeq} (in the hyperk\"ahler case) by a formal substitution in which $\K$ becomes a scalar field 
and the term $(\nabla_{A'}\K)(\nabla^{A'}\K)$, which is a commutator of matrices, 
is replaced by certain Poisson bracket, see e.g. \cite{Monteiro, Chacon}.
(We are not aware of a similar formal procedure for the hyper-heavenly case, but it might be interesting 
to analyse this possibility.)
This analogy between the equation satisfied by a SDYM field and the one satisfied by a 
SD gravitational field is an instance of what is known as the double copy,
a conjectured relationship between scattering amplitudes for gauge and gravity theories that is of much 
current interest. In this regard, we note the following: 
\begin{itemize}
\item The similarity between the equations for Yang-Mills fields and for gravity in this note originates simply
in the fact that an Einstein metric (not necessarily self-dual) can be regarded as a particular case of a SDYM system, 
a result that follows from the work of Atiyah, Hitchin and Singer in \cite{Atiyah}. We will give an explicit version of this below.
\item For the gravitational case the internal indices in \eqref{Kmatrixeq} (the matrix indices of $\K$)
can be converted to spinor indices: see eq. \eqref{SpinorConnection}.
In other words, internal indices are converted to ``external'' indices, since spinors are indeed related 
to directions in spacetime. 
The internal two-dimensional Yang-Mills space becomes the ``external'' spinor space. 
This is, of course, independent of any field equation or any special background geometric structure. 
It is ultimately associated to the fact that (in any signature) the Riemann tensor 
is an element of $\Lambda^{2}\otimes\Lambda^{2}$ (the ``square'' of $\Lambda^{2}$), see the section below.
\end{itemize}

\subsubsection*{The Einstein equations}

So far the vector bundle $E$ and the connection ${\rm D}$ have been arbitrary, representing a generic Yang-Mills system. 
Following the result of Atiyah-Hitchin-Singer in \cite[Proposition 2.2]{Atiyah}, 
we can interpret the Einstein equations as a particular case of a SDYM field
with some additional structures on the vector bundle. 
The essential idea is as follows\footnote{We focus on Riemann signature, but the same result can be proved 
in Lorentz and split signature as well.}.
Given an orientable Riemannian manifold $(M,g)$, we can construct (at least locally) 
the spinor bundles $\mathbb{S}$, $\mathbb{S}'$. 
For the argument it is sufficient to consider only one of them, say $\mathbb{S}'$.
This is a rank-2 vector bundle over $M$, equipped with a symplectic form $\epsilon'$. 
The transition functions then take values in the symplectic group ${\rm Sp}(1)\cong{\rm SU}(2)$.
The Levi-Civita connection of $g$ induces the spin connection 
$\nabla:\Gamma(\mathbb{S}')\to\Gamma(T^{*}M\otimes\mathbb{S}')$, which is compatible 
with $\epsilon'$ in the sense that $\nabla\epsilon'=0$.
One can compute the curvature of this connection, say $\mathcal{R}^{+}\in\Gamma(\Lambda^{2}\otimes{\rm End}(\mathbb{S'}))$, 
and it turns out that this is SD if and only if the metric is Einstein.
For convenience, below we give an explicit proof of this.

Let $\mathcal{R}$ denote the Riemann curvature tensor of $g$ with 
all its indices in the lower position. From the well-known algebraic symmetries of $\mathcal{R}$, we have 
$\mathcal{R}\in\Lambda^{2}\otimes\Lambda^{2}$. (For simplicity we omit the symbol `$\Gamma$' corresponding to 
spaces of sections.)
Recalling the decomposition $\Lambda^{2}=\Lambda^{2}_{+}\oplus\Lambda^{2}_{-}$, 
and applying this to the second factor in $\Lambda^{2}\otimes\Lambda^{2}$, 
one gets $\mathcal{R}=\mathcal{R}^{+}+\mathcal{R}^{-}$, where $\mathcal{R}^{\pm}\in\Lambda^{2}\otimes\Lambda^{2}_{\pm}$.
Decomposing now the first $\Lambda^{2}$ factor: $\mathcal{R}^{+}=\tilde{X}+\tilde{\Phi}$ and $\mathcal{R}^{-}=X+\Phi$, 
where $\tilde{X}\in\Lambda^{2}_{+}\otimes\Lambda^{2}_{+}$, $\tilde{\Phi}\in\Lambda^{2}_{-}\otimes\Lambda^{2}_{+}$, 
$\Phi\in\Lambda^{2}_{+}\otimes\Lambda^{2}_{-}$ and $X\in\Lambda^{2}_{-}\otimes\Lambda^{2}_{-}$.
The tensor $\tilde{X}$ contains the SD Weyl tensor and the curvature scalar, and $\tilde{\Phi}$ corresponds 
to the trace-free Ricci tensor. Similarly, $X$ contains the ASD Weyl tensor and the curvature scalar, 
and $\Phi$ also corresponds to the trace-free Ricci tensor. In practice, $\tilde{\Phi}$ and $\Phi$ contain the same 
information\footnote{This is a consequence of the extra symmetry $R_{abcd}=R_{cdab}$ of the Riemann tensor.}; 
here we distinguish them because of their duality properties. 
The metric $g$ is Einstein if and only if $\tilde{\Phi}=0$. 
On the other hand, we can dualize the first form indices in $\mathcal{R}$, this gives 
${}^{*}\mathcal{R}^{+}=\tilde{X}-\tilde{\Phi}$ and ${}^{*}\mathcal{R}^{-}=-X+\Phi$.
Therefore, we see that ${}^{*}\mathcal{R}^{+}=\mathcal{R}^{+}$ if and only if $\tilde{\Phi}=0$, 
that is, the curvature tensor $\mathcal{R}^{+}$ is SD if and only if the metric is Einstein. 
The fact that $\mathcal{R}^{+}$ is the curvature of the connection on the spin bundle $\mathbb{S}'$ 
follows from the isomorphism between $\Lambda^{2}_{+}$ and (two copies of) $\mathbb{S}'$.

\medskip
In abstract indices, the curvature decomposition above is 
$X_{abcd}=X_{ABCD}\epsilon_{A'B'}\epsilon_{C'D'}$, $\tilde{X}_{abcd}=\tilde{X}_{A'B'C'D'}\epsilon_{AB}\epsilon_{CD}$, 
$\Phi_{abcd}=\Phi_{A'B'CD}\epsilon_{AB}\epsilon_{C'D'}$ and $\tilde{\Phi}_{abcd}=\Phi_{ABC'D'}\epsilon_{A'B'}\epsilon_{CD}$, 
where all spinors are defined as in \cite{PR1}. 
The tensor field $\mathcal{R}^{+}$ is $\mathcal{R}^{+}_{abcd}=R^{+}_{abC'D'}\epsilon_{CD}$, where
$R^{+}_{abC'D'}=\tilde{X}_{A'B'C'D'}\epsilon_{AB}+\Phi_{ABC'D'}\epsilon_{A'B'}$.
With this notation, the Einstein equations are equivalent to the following self-duality condition:
\begin{align}
  {}^{*}R^{+}_{abC'D'} = R^{+}_{abC'D'}
  \quad \Leftrightarrow \quad \Phi_{ABC'D'}=0.
\end{align}
The first two SDYM equations in \eqref{SDYMcurvedspace}, $F_{zw}=0$ and $F_{\tilde{z}\tilde{w}}=0$, are 
now $\iota^{A}\iota^{B}\Phi_{ABC'D'}=0$ and $o^{A}o^{B}\Phi_{ABC'D'}=0$ respectively.
From the previous general discussion, they imply the existence of local holomorphic and anti-holomorphic 
trivializations of $\mathbb{S}'$.
In the terminology and notation of this subsection, this means that there are two (not necessarily normalized) 
primed spin frames, $\{e^{A'}_{\bf i}\}$ and $\{\tilde{e}^{A'}_{\bf i}\}$, with ${\bf i}=0,1$, such that 
$\nabla_{A'}e^{B'}_{\bf i} = 0$, $\tilde{\nabla}_{A'}\tilde{e}^{B'}_{\bf i} = 0$.
Yang's matrix is the $2\times2$ matrix defined by $e^{A'}_{\bf j} = \J^{\bf i}{}_{\bf j}\tilde{e}^{A'}_{\bf i}$, 
and the generalized Yang equation \eqref{YangEqCurvedSpace} is 
\begin{equation}
 \nabla^{A'}[r(\J^{-1})^{\bf i}{}_{\bf k}\tilde{\nabla}_{A'}\J^{\bf k}{}_{\bf j}] = 0. \label{YangEqEinstein}
\end{equation}
The $\K$-matrix equation \eqref{Kmatrixeq} has the same index structure.

We note that in the hyperk\"ahler case, a calculation shows that Yang's matrix $\J$ is essentially 
the Hermitian matrix $g_{\mu\tilde{\nu}}$, see the appendix \ref{App:2HE}.

\medskip
Before studying B\"acklund transformations for the general SDYM equations, 
let us make an additional comment on the relation between Einstein metrics and SDYM systems. 
We just argued that, if we start with a Riemannian manifold whose metric tensor is Einstein,
then we can view this as a SDYM system for the spin connection on primed spinors.
It is worth remarking, however, that the `converse' construction is more subtle. 
Namely, suppose we start with a Riemmanian metric and a rank-2 vector bundle equipped 
with a SD connection, can we construct an Einstein metric from such a connection?
In other words, in which situations, or with which additional structures, 
can we interpret the SD connection as the spin connection on the primed spinor bundle of an Einstein metric?
Unfortunately, the answer to this question is not clear to us;
but it would be useful for the application of proposition \ref{Prop:hiddensymmetries} below
to the generation of new Einstein metrics.

\section{B\"acklund transformations}

Having obtained a generalization \eqref{YangEqCurvedSpace} of the flat space Yang's equation \eqref{YangEqFlatspace}, 
we now wish to investigate if the hidden symmetry structure of \eqref{YangEqFlatspace} can be extended to the curved case. 
We follow the procedure of Mason and Woodhouse \cite[Section 4.6]{MW} and write $\J$ in the form
\begin{equation}
 \J = 
 \left(  \begin{matrix}
 \mathsf{A}^{-1} - \tilde{\mathsf{B}}\tilde{\mathsf{A}}\mathsf{B} & -\tilde{\mathsf{B}}\tilde{\mathsf{A}} \\
 \tilde{\mathsf{A}}\mathsf{B} & \tilde{\mathsf{A}}
 \end{matrix} \right) 
 = 
  \left(  \begin{matrix}
 1  & \tilde{\mathsf{B}} \\
 0 & \tilde{\mathsf{A}}^{-1}
 \end{matrix} \right)^{-1} 
  \left(  \begin{matrix}
 \mathsf{A}^{-1}  & 0 \\
 \mathsf{B} & 1
 \end{matrix} \right)
 \label{decompositionJ}
\end{equation}
where $\mathsf{A}$ is a $k\times k$ matrix, $\tilde{\mathsf{A}}$ is $\tilde{k}\times\tilde{k}$ (with $k+\tilde{k}=n$), 
$\mathsf{B}$ is $\tilde{k}\times k$, and $\tilde{\mathsf{B}}$ is $k\times\tilde{k}$.
In flat space, the decomposition as a product $\J\equiv\tilde{h}^{-1}h$ in the second equality in \eqref{decompositionJ}
allows to relate $\J$ to a patching matrix in twistor space. (Here we find the decomposition useful for computing $\J^{-1}$.)
The generalized Yang equation \eqref{YangEqCurvedSpace} becomes then the following system of four equations:
\begin{subequations}
\begin{align}
 & \nabla^{A'}[r\sA(\tilde{\nabla}_{A'}\tilde{\sB})\tilde{\sA}] = 0, \label{yangBacklund1} \\
 & \nabla^{A'}[r\sB\sA(\tilde{\nabla}_{A'}\tilde{\sB})\tilde{\sA} + r\tilde{\sA}^{-1}\tilde{\nabla}_{A'}\tilde{\sA}] = 0, \label{yangBacklund2} \\
 & \nabla^{A'}[r\sA(\tilde{\nabla}_{A'}\tilde{\sB})\tilde{\sA}\sB + r(\tilde{\nabla}_{A'}\sA)\sA^{-1} ] = 0, \label{yangBacklund3} \\
 & \nabla^{A'}[r\sB\sA(\tilde{\nabla}_{A'}\tilde{\sB})\tilde{\sA}\sB 
 + r(\tilde{\sA}^{-1}\tilde{\nabla}_{A'}\tilde{\sA})\sB + r\sB(\tilde{\nabla}_{A'}\sA)\sA^{-1} + r\tilde{\nabla}_{A'}\sB] = 0. 
 \label{yangBacklund4}
\end{align}
\end{subequations}
Using the first three equations, together with the identity $\nabla^{A'}(r\tilde{\nabla}_{A'}\psi)=\tilde{\nabla}_{A'}(r\nabla^{A'}\psi)$
valid for any matrix-valued function $\psi$, the fourth equation can be arranged conveniently so that we get
$\tilde{\nabla}^{A'}[r\tilde{\sA}(\nabla_{A'}\sB)\sA] = 0$.

\smallskip
Now we need to apply an integration procedure\footnote{The following is a simple version of the construction 
presented in \cite[Section 2.5]{Araneda}, which is for connections on more general vector bundles.}
associated to the operators \eqref{partialconnections}:
if a matrix-valued spinor $a_{A'}$ satisfies $\tilde{\nabla}^{A'}a_{A'}=0$, 
then there exists, locally, a matrix-valued function $a$ such that $a_{A'}=\tilde{\nabla}_{A'}a$; and analogously for 
an equation of the form $\nabla^{A'}b_{A'}=0$ (since $\nabla_{A'}$ is the complex conjugate of $\tilde{\nabla}_{A'}$).
This can be seen as follows. 
First write $0=\epsilon^{A'B'}\tilde{\nabla}_{A'}a_{B'}$ and express this in terms of the 
coordinate vector fields $\partial_{\tilde{z}}, \partial_{\tilde{w}}$, using
the identity $\epsilon^{A'B'} = \tilde{N}^{-1}(\tilde{Z}^{A'}\tilde{W}^{B'} - \tilde{W}^{A'}\tilde{Z}^{B'})$, 
where $\tilde{Z}^{A'},\tilde{W}^{A'}$ were defined in \eqref{CoordinateVF1}, \eqref{CoordinateVF2}.
Then use the fact that $[\partial_{\tilde{z}}, \partial_{\tilde{w}}]=0$, 
which implies $\nabla_{\partial_{\tilde{z}}}\tilde{W}^{B'}-\nabla_{\partial_{\tilde{w}}}\tilde{Z}^{B'}=0$.
Finally, the equation $\partial_{\tilde{z}}a_{\tilde{w}}-\partial_{\tilde{w}}a_{\tilde{z}}=0$
implies that there exists, locally, a matrix-valued function $a$ such that 
$a_{\tilde{w}} = \partial_{\tilde{w}}a$, $a_{\tilde{z}}=\partial_{\tilde{z}}a$, or equivalently $a_{A'}=\tilde{\nabla}_{A'}a$.

\smallskip
Applying the above procedure to equations \eqref{yangBacklund1} and $\tilde{\nabla}^{A'}[r\tilde{\sA}(\nabla_{A'}\sB)\sA] = 0$,
we deduce that there are two matrices $\sB'$, $\tilde{\sB}'$ such that
\begin{subequations}
\begin{align}
 & r\sA(\tilde{\nabla}_{A'}\tilde{\sB})\tilde{\sA} = \nabla_{A'}\sB', \label{B'} \\
 & r\tilde{\sA}(\nabla_{A'}\sB)\sA = \tilde{\nabla}_{A'}\tilde{\sB}'. \label{tB'}
\end{align}
\end{subequations}
We note that the matrices $\sB'$, $\tilde{\sB}'$ are not uniquely defined: we have the freedom $\sB'\to \sB'+\sC'$, 
$\tilde{\sB}' \to \tilde{\sB}' + \tilde{\sC'}$ where $\nabla_{A'}\sC'=0$, $\tilde{\nabla}_{A'}\tilde{\sC}'=0$
(so $\sC$ and $\sC'$ are a sort of `integration constants').
Now define the discrete (`B\"acklund') transformation
\begin{equation}
 (\sA,\tilde{\sA},\sB,\tilde{\sB},k,\tilde{k}) \to (\sA',\tilde{\sA}',\sB',\tilde{\sB}',\tilde{k},k) \label{BacklundT}
\end{equation}
where for the moment $\sA'$ and $\tilde{\sA}'$ are undetermined, and let $\J'$ be 
a new matrix defined as \eqref{decompositionJ} but with the primed matrices instead.
Yang's equation \eqref{YangEqCurvedSpace} for $\J'$ is then the system 
\eqref{yangBacklund1}--\eqref{yangBacklund4} for the primed matrices. 
The primed equation \eqref{yangBacklund1} gives $\nabla^{A'}[r^{2}\sA'\tilde{\sA}(\nabla_{A'}\sB)\sA\tilde{\sA}']$, 
which suggests us to set $\sA'\equiv r^{-1}\tilde{\sA}^{-1}$ and $\tilde{\sA}'\equiv r^{-1}\sA^{-1}$. 
(The ordinary B\"acklund transformation in flat space is defined by \eqref{BacklundT} but 
with $\sA'=\tilde{\sA}^{-1}$ and $\tilde{\sA}'=\sA^{-1}$, \cite[Section 4.6]{MW}.)
Then we find
\begin{align*}
 & \nabla^{A'}[r\sA'(\tilde{\nabla}_{A'}\tilde{\sB}')\tilde{\sA}'] = 0,  \\
 & \nabla^{A'}[r\sB'\sA'(\tilde{\nabla}_{A'}\tilde{\sB}')\tilde{\sA}' + r\tilde{\sA}'^{-1}\tilde{\nabla}_{A'}\tilde{\sA}'] = -\nabla^{A'}\tilde{\nabla}_{A'}r,  \\
 & \nabla^{A'}[r\sA'(\tilde{\nabla}_{A'}\tilde{\sB}')\tilde{\sA}'\sB' + r(\tilde{\nabla}_{A'}\sA')\sA'^{-1} ] = -\nabla^{A'}\tilde{\nabla}_{A'}r,  \\
 & \nabla^{A'}[r\sB'\sA'(\tilde{\nabla}_{A'}\tilde{\sB}')\tilde{\sA}'\sB' 
 + r(\tilde{\sA}'^{-1}\tilde{\nabla}_{A'}\tilde{\sA}')\sB' + r\sB'(\tilde{\nabla}_{A'}\sA')\sA'^{-1} 
 + r\tilde{\nabla}_{A'}\sB'] = -2(\nabla^{A'}\tilde{\nabla}_{A'}r) \sB'
\end{align*}
Thus, $\J'$ will satisfy Yang's equation \eqref{YangEqCurvedSpace} if and only if the scalar field $r$ (defined in \eqref{r})
satisfies $\nabla^{A'}\tilde{\nabla}_{A'}r=0$.
However, a lengthy calculation gives $\nabla^{A'}\tilde{\nabla}_{A'}r=-3r\chi^{-1}\Psi_{2}$, where 
$\Psi_{2}=\Psi_{ABCD}o^{A}o^{B}\iota^{C}\iota^{D}$;
therefore: 
\begin{align}
 \nabla^{A'}\tilde{\nabla}_{A'}r=0 \quad \Leftrightarrow \quad \Psi_{2}\equiv 0.
\end{align}
Given that $\Psi_{2}$ is the only non-trivial component of $\Psi_{ABCD}$ (recall that the conformal K\"ahler assumption 
implies that $\Psi_{ABCD}$ is type D),
we see that the new $\J'$ will satisfy Yang's equation \eqref{YangEqCurvedSpace} if and only if $\Psi_{ABCD}\equiv 0$.

The above discussion then produces the following result:

\begin{proposition}\label{Prop:hiddensymmetries}
Let $(M,g,J)$ be a conformally-K\"ahler, conformally half-flat Riemannian manifold\footnote{This is equivalent 
to saying that $g$ is conformal to a scalar-flat K\"ahler metric, see e.g. \cite[Proposition 13.4.9]{MW}.}. 
Let $\J$ be the Yang's matrix of a SDYM field, expressed in the form \eqref{decompositionJ}, 
and let $\J'$ be the matrix resulting from the B\"acklund transformation \eqref{BacklundT}. 
Then $\J'$ is the Yang's matrix of a new solution to the SDYM equations.
Combining this transformation with conjugation by constant matrices, we get a continuous symmetry group 
for the SDYM equations. 
\end{proposition}

One may argue, of course, that there might exist a transformation different from \eqref{BacklundT} 
that produces a new SDYM field without having to restrict to $\Psi_{ABCD}\equiv 0$. 
Unfortunately, after exploring several possibilities, we were not able to avoid the restriction $\Psi_{ABCD}=0$. 

\smallskip
As we saw, this solution generating technique can also employ certain Einstein metrics as seeds.
More precisely, suppose we have a conformally half-flat Einstein space,
which is by definition a quaternionic-K\"ahler metric. 
In spinors, this means that $\Psi_{ABCD}=0$, $\Phi_{ABC'D'}=0$, while the scalar curvature $\Lambda$ 
and the SD Weyl spinor $\tilde{\Psi}_{A'B'C'D'}$ are arbitrary. 
Then we can construct a $2\times2$ Yang's matrix $\J$, which will satisfy eq. \eqref{YangEqEinstein}.
The B\"acklund transformation \eqref{BacklundT} generates a new matrix $\J'$ which 
also satisfies \eqref{YangEqEinstein}. 
Then the gauge potential \eqref{gaugePotentialSpinors} gives a new solution to 
the SDYM equations, with gauge group ${\rm GL}(2,\mathbb{C})$.
We note, however, that from our discussion of the Einstein equations in section \ref{Sec:curvedspaces}, 
it is a priori not possible to conclude that the new ${\rm GL}(2,\mathbb{C})$-SDYM field 
gives a new solution to the Einstein equations.

\appendix

\section{The $\K$-matrix and the second heavenly equation}\label{App:2HE}

In this appendix we consider a (4D) hyperk\"ahler manifold $(M,g_{ab})$, and we want to show that \eqref{Kmatrixeq}
reduces to Pleba\'nski's second heavenly equation.
A geometry is hyperk\"ahler if it is Ricci-flat and self-dual. In spinors, this means 
that $\Psi_{ABCD}=\Phi_{ABC'D'}=\Lambda=0$ and $\tilde{\Psi}_{A'B'C'D'}$ is arbitrary.
An equivalent characterization is to say that there is a parallel spinor: $\nabla_{AA'}o_{B}=0$. 
The complex conjugate $\iota_{B}=o^{\dagger}_{B}$ is then also parallel, and together they 
give a parallel unprimed spin frame. 
With respect to the notation in section \ref{Sec:curvedspaces}, we can simply take $o_{A}\iota^{A}=\chi=1$ and $\phi=1$.
The complex structure \eqref{ComplexStructureSpinors} is of course integrable, with holormorphic coordinates $(z,w)$.

\smallskip
In order to arrive to eq. \eqref{Kmatrixeq}, we need two spin frames $e^{A'}_{\bf i}$ and $\tilde{e}^{A'}_{\bf i}$ 
parallel under $\nabla_{A'}$ and $\tilde{\nabla}_{A'}$, respectively.
Notice that we actually need only one of them, say $e^{A'}_{\bf i}$, since then its complex conjugate 
will be parallel under $\tilde{\nabla}_{A'}$.
In the hyperk\"ahler case we can obtain this frame from the definition of complex coordinates, as follows.
Equation \eqref{Coordinates1forms} for the complex structure $J_{0}$ gives 
${\rm d}z=o_{A}\alpha_{A'}{\rm d}x^{AA'}$ and ${\rm d}w=o_{A}\beta_{A'}{\rm d}x^{AA'}$, 
where $\alpha_{A'}\beta^{A'}=k\neq 0$.
From the fact that ${\rm d}^{2}z=0={\rm d}^{2}w$, a straightforward calculation
gives $\tilde{\nabla}_{A'}\alpha_{B'}=0$, $\tilde{\nabla}_{A'}\beta_{B'}=0$, and also 
$\nabla_{A'}\alpha^{A'}=0=\nabla_{A'}\beta^{A'}$.
This implies that $\tilde{\nabla}_{A'}k=0$, or in other words $k=k(z,w)$.
Using the freedom in the definition of the coordinates $(z,w)$, namely $z\to z'=z'(z,w)$, $w\to w'=w'(z,w)$,
we can then choose coordinates so that $k\to k'=1$. 
Omitting the prime, we have $e^{A'}_{0}:=\alpha^{\dagger A'}$, $e^{A'}_{1}:=\beta^{\dagger A'}$, 
$e_{0A'}e^{A'}_{1}=1$, $\nabla_{A'}e^{B'}_{\bf i}=0$ and $\tilde{\nabla}_{A'}e^{A'}_{\bf i}=0$.
In addition, we take $\tilde{e}^{A'}_{\bf i} = e^{\dagger A'}_{\bf i}$.

\smallskip
We note that, using the definition of Yang's matrix, i.e. $e^{A'}_{\bf j} = \J^{\bf i}{}_{\bf j}\tilde{e}^{A'}_{\bf i}$, 
together with $g_{z\tilde{z}}=g(\partial_{z},\partial_{\tilde{z}})$, etc., 
a short computation gives the components $\J^{0}{}_{0}=-g_{z\tilde{w}}$, $\J^{1}{}_{0}=-g_{w\tilde{w}}$, 
$\J^{0}{}_{1}=g_{z\tilde{z}}$, $\J^{1}{}_{1} = g_{w\tilde{z}}$, so $\J$ is essentially given by the Hermitian matrix $g_{\mu\tilde{\nu}}$.

\smallskip
Let $\theta^{\bf i}_{A'}$ be such that $\theta^{\bf i}_{A'}e^{A'}_{\bf j}=\delta^{\bf i}{}_{\bf j}$. 
(Note that if we define $e^{\bf i}_{A'}:=\epsilon^{\bf ij}\epsilon_{B'A'}e^{B'}_{\bf j}$, then $\theta^{\bf i}_{A'}=-e^{\bf i}_{A'}$.)
The connection 1-form is given by $\nabla_{AA'}e^{B'}_{\bf i}=A_{AA'}{}^{\bf j}{}_{\bf i}e^{B'}_{\bf j}$.
The internal indices ${\bf i,j}$ can be converted to ``external'' by defining 
\begin{equation}
 A_{AA'B'C'} \equiv \theta^{\bf i}_{C'}\nabla_{AA'}e_{{\bf i}B'}. \label{SpinorConnection}
\end{equation}
This object satisfies $A_{AA'B'C'}=A_{AA'(B'C')}$ and $\iota^{A}A_{AA'B'C'}=0$.
In addition, the fact that $\tilde{\nabla}_{A'}e^{A'}_{\bf i}=0$ implies $A_{AA'B'C'}=A_{A(A'B')C'}$.
Therefore $A_{AA'B'C'}$ is totally symmetric in the primed indices:
\begin{equation}
 A_{AA'B'C'} = A_{A(A'B'C')}. \label{totallysymmetric}
\end{equation}

Now, the Ricci spinor is 
\begin{equation*}
 \Phi_{ABC'D'} = \nabla_{A'(A}A^{A'}{}_{B)C'D'} + A_{A'(A}{}^{E}{}_{|C'|}A^{A'}{}_{B)D'E'}.
\end{equation*}
The equation $\iota^{A}\iota^{B}\Phi_{ABC'D'}=0$ is automatic since $\iota^{A}A_{AA'B'C'}=0$.
From $o^{A}\iota^{B}\Phi_{ABC'D'}=0$ we immediately get $\nabla^{A'}(o^{A}A_{AA'B'C'})=0$, 
which, in view of \eqref{totallysymmetric}, implies that there is a scalar function $\Theta$ such that 
\begin{equation}
 A_{AA'B'C'} = -\iota_{A}\nabla_{A'}\nabla_{B'}\nabla_{C'}\Theta. \label{ScalarPotential}
\end{equation}
The SD Weyl spinor is $\tilde{\Psi}_{A'B'C'D'}=\nabla_{A(A'}A^{A}{}_{B'C'D')}=-\nabla_{A'}\nabla_{B'}\nabla_{C'}\nabla_{D'}\Theta$.
Replacing these identities in the last equation $o^{A}o^{B}\Phi_{ABC'D'}=0$, after rearranging terms one arrives at
\begin{equation*}
 \nabla_{C'}\nabla_{D'} \left[ \Box\Theta + (\nabla_{A'}\nabla_{B'}\Theta)(\nabla^{A'}\nabla^{B'}\Theta) \right] = 0
\end{equation*}
which is the same as the $\mathcal{K}$-matrix equation \eqref{Kmatrixeq}. Equivalently:
\begin{equation}
 \Box\Theta + (\nabla_{A'}\nabla_{B'}\Theta)(\nabla^{A'}\nabla^{B'}\Theta) = C \label{2HE}
\end{equation}
where $C$ is a function such that $\nabla_{C'}\nabla_{D'}C=0$. 
This can be absorbed in the definition of $\Theta$ so that we get zero in the RHS of \eqref{2HE};
the result is Pleba\'nski's second heavenly equation.

\end{document}